\documentstyle[12pt]{article}
\title{ The CUSUM  test for detecting structural changes in strong mixing processes}
\author{Fatemeh Azizzadeh\thanks{School of Mathematics and Computer Science,
Amirkabir University of Technology, Tehran, Iran}, S.
Rezakhah\thanks{Faculty of Mathematics and Computer Science,
Amirkabir University of Technology, Tehran, Iran.
 Email: rezakhah@aut.ac.ir, fatemeh_aziz@aut.ac.ir}}

\headheight\baselineskip \headsep2\baselineskip \textwidth200mm
\advance\textwidth -1.6in \textheight275mm
\advance\textheight-1.7in \footskip2\baselineskip

\begin{document}
\pagestyle{plain}
\title{ The CUSUM  test for detecting structural changes in strong mixing processes}
\author{F. Azizzadeh, S. Rezakhah \thanks{ Faculty
of Mathematics and Computer Science, Amirkabir University of
Technology, Tehran, Iran. Email: rezakhah@aut.ac.ir, fatemeh\_aziz@aut.ac.ir }}
\maketitle
\renewcommand{\theequation}{\arabic{section}.\arabic{equation}}
\begin{abstract}
Strong mixing property holds for a broad class of linear and nonlinear
time series models such as ARMA and GARCH models.
In this article we study correlation structure of
strong mixing sequences, and some asymptotic properties are presented.
We also present a new method for detecting change point in correlation structure of strong mixing sequences, and
 present a nonparametric CUSUM test statistic for this. Asymptotic
 consistency of this test statistics is shown. This method is applied to simulated data
  of some linear and nonlinear models and power of the test
 is evaluated. For linear models, it is shown that
  this method have a better performance in compare to Berkes et al.(2009).
 \end{abstract}
 MSC: Primary 62M10, 60F17, Secondary 62G20, 62G10.\\
 \noindent Key words and Phrases: structural changes, Strong mixing,
 Functional central limit theorem, CUSUM test, Brownian bridge.\\
\setcounter{equation}{0}
\newcommand{\be}{\begin{equation}} \newcommand{\ee}{\end{equation}}

\section{Introduction}
\hspace{.2in}Change point detection in a sequence of random
variables was first proposed by Page(1954).
This study started by detecting changes in
the mean of a sequence of independent random variables and then
extended to dependent sequences. Change point detection is widely
 used in various fields such as quality control, economics,
finance and medicine. Review
of earlier works can be found in Cs\"{o}rg\"{o} and
Horv\'{a}th(1988), Brodsky and Darkhovsky (1993) and
Cs\"{o}rg\"{o} and Horv\'{a}th(1997). \par
\hspace{.2in}Among different methods for change point detection,
 the CUSUM test proposed by Page(1954), for mean change detection,
  is widely used for its simplicity.
Incl\'{a}n and Tiao(1994) proposed a CUSUM of squares test for testing a variance change in i.i.d. normal random variables.
Lee and Park(2001) extended the CUSUM test of squares test of Incl\'{a}n and
Tiao(1994) for linear processes. Lee et al.(2003) studied
change of parameters in a random coefficient AR(1)
model, thus detecting changes in the auto-covariances of a linear
process
Galeano and Pena(2007) studied changes in variance and correlation
structure of the multivariate time series. Zhou and Liu(2009) used
a weighted CUSUM statistic for mean change detection in infinite variance AR(p)
process.  Berkes et al.(2009) considered a CUSUM test to detect changes
in the mean and in the covariance structure of a linear process.
Recently Qin et al.(2010) studied mean change detection in $\alpha$-mixing
processes.\par
\hspace{.2in} In this article we study change in the correlation structure of strong
mixing sequences. Let $\{X_t: t\geq 1\}$ be a stationary process.
As a measure of dependence we use Rosenblatt's $\alpha$-mixing
coefficient as
$$
\alpha_X(n, j)=\sup_{A,B}|P(A\cap B)-P(A)P(B)|, \eqno(1.1)
$$
where A and B are in the $\sigma$-fields
${\mathbf{M}}_{-\infty}^{n}(X)=\sigma\{X_t; t<n\}$ and
${\mathbf{M}}_{n+j}^{\infty}(X)=\sigma\{X_t; t>n+j\}$ respectively.
The sequence $\{X_t\}$ is said to be $\alpha$-mixing or strong
mixing (SM) if
$$
\alpha_X(j)= \sup_n\alpha_X(n, j) \rightarrow 0 \ \ as\
j\rightarrow\infty.  \eqno(1.2)
$$
\hspace{.2in}Strong mixing processes are asymptotically independent.
Strong mixing property holds for a large class of linear
and nonlinear stationary time series such as ARMA and GARCH models,
 m-dependent processes, broad class of
Gaussian processes and ergodic Markov processes( Bosq 1996, and Bradley 2005). \par
\hspace{.2in} Ibragimov(1962) showed some results for stationary strong
mixing sequences and proved central limit theorem for strict stationary
SM processes. Davydov(1968) obtained some
moment inequalities and Rio(1993) presented some covariance
inequalities and bounds on the variance of partial sums of
SM processes. Herrndorf(1985), Doukhan et al.(1994), and
Merlevede and Peligard(2000) studied functional central
limit theorem on SM processes. Romano
and Thombs(1996) used central limit theorem, established by
Ibragimov(1962), to show that sample auto-covariances of strictly
stationary SM sequences converge in distribution to normal
distribution.\par
\hspace{.2in}By using functional central limit theorem for SM
sequences, we propose a new test statistic for detecting
changes in correlation structure of stationary strong
mixing processes.
This test is a nonparametric one and does not depend on any
assumptions about the underlying distribution or model.
 The rest of this paper is organized as follows:
In section 2  a nonparametric test statistic for detection of change points, in broad class of linear and non linear process, is constructed
and its asymptotic properties, under no change null hypothesis, are studied.
Also the consistency of this test statistic is shown
 in section 2. Section 3 is devoted to the simulation results on
 different linear and nonlinear models. In this section
  the method of Berkes et al.(2009) has been compared
  with the method of this paper  for some some linear models.
 By simulation we show that this test statistic have
 a better performance and is more powerful in many cases.
\section{Main results}
\hspace{.2in}In this section we present some preliminary results which will be used later in this paper.
 We  present functional central limit
theorem for sample auto-covariances of SM processes. We also introduce
 a new test statistic for detecting changes in the correlation
 structure of stationary SM processes,
 which we call CUSUM strong mixing(CSSM).
 Finally we show consistency and
asymptotic convergence of this empirical CSSM test statistic.\par
\hspace{.2in}Let $\{ X_n\}$ be a sequence of random variables on some
probability space $(\Omega, \mathcal{F}, P)$, satisfying
$$
E(X_n)=0, \ \ \ E(X_n^2)<\infty \ \ \ \ \mbox{for} \ \ n>0. \eqno(2.1)
$$
\hspace{.2in}Let $S_n=X_1 + ... + X_n$ for $n>0$.
Consider the Skorokhod space $D\equiv D[0, 1]$ of all functions on $[0, 1]$
which are right continuous with left limit. Let
 $W_n(t):\Omega \rightarrow D$ to be a random function as
$$
W_n(t)=\frac{S_{[nt]}}{\sigma \sqrt{n}}\ \ \ \ \mbox{for} \ \ t\in[0, 1],
\ n>0.
$$
If $W_n(t)$ is weakly convergent to a standard Brownian motion $W(t)$,
then $X_n$ is said to satisfy the functional central limit
 theorem or strong invariance principle (Billingsley, 1999).
Herrndorf(1985) proved functional central limit theorem for strong
mixing sequences without stationarity assumption but
 assumed convergence of the variance of partial sums. \par
\hspace{.2in}For $X_1, X_2, ...$ as a sequence of zero-mean stationary process
 the sample auto-covariances
$\hat{\gamma}_n(h)$, $h=0, ..., n$, are defined as:
$$
\hat{\gamma}_n(h)= \hat{\gamma}_n(-h)=
\frac{1}{n}\sum_{i=1}^{n-h}X_iX_{i+h}. \eqno(2.2)
$$
Asymptotic covariance of sample autocovariances is known as
Bartlett's estimator and is defined as:
$$
c_{h k}= \lim_{n\rightarrow \infty}nCov(\hat{\gamma}_n(h),
\hat{\gamma}_n(k)). \eqno(2.3)
$$
Let
$$
g^*_n(t):=\frac{[nt]}{\sqrt{n}}C^{-1/2} (
\hat{\gamma}_{[nt]}(0)-\hat{\gamma}_n(0),\,
 \hat{\gamma}_{[nt]}(1)-\hat{\gamma}_n(1), ...,
\, \hat{\gamma}_{[nt]}(L)-\hat{\gamma}_n(L))^T, \eqno(2.4)
 $$
where  $C=[c_{hk}]_{h, k=1}^{L+1}$
is the covariance matrix whose entries are defined by (2.3). Now we have the following result.
\\
\emph{\noindent
\textbf{Theorem 1}: Let $\{X_n\}_{n=1}^\infty$ be a stationary
strong mixing process which satisfies:
\begin{itemize}
  \item i) $\sup_i E|X_i|^{4+2\delta}<\infty$
  \item ii) $\sum_{k=1}^{\infty}\alpha_X(k)^{\delta/(2+\delta)}<\infty$
   , for some $\delta\in (0, \infty)$
\end{itemize}
in which $\alpha_X(.)$ is the mixing coefficient, defined by (1.2). Then
$$
 g^*_n(t)^T . g^*_n(t) \Rightarrow \sum_{j=0}^{L}(W^0_j(t))^2 \ \ \ \
  \mbox{as} \ n \rightarrow \infty, \eqno(2.5)
$$
where $\Rightarrow$ denotes convergence in distribution, $W^0_j(.)$
are independent Brownian bridge, for $0\leq j\leq L$.}
\\

Before proceeding to the proof of this theorem we present some lemmas,
which are necessary for our proof. \\

\emph{\noindent
\textbf{Lemma 1}: (Davydov, 1968) Let the process $\{X_t\}$
be strong mixing, and random variables
$\xi$ and $\eta$ be measurable with respect to
${\mathbf{M}}_{-\infty}^{n}(X)$ and ${\mathbf{M}}_{n+j}^{\infty}(X)$,
 introduced by (1.1), respectively.
Moreover if the moments $E|\xi|^p$ and $E|\eta|^q$
exist for $p,q>0$ where $\frac{1}{p}+\frac{1}{q}<1$, then
$$
|E\xi\eta- E\xi E\eta|\leq
C[E|\xi|^p]^{1/p}[E|\eta|^q]^{1/q}[\alpha(n)]^{1-1/p-1/q}.
$$
}
\emph{\noindent
\textbf{Lemma 2}: Let $\{X_n\}_{n=1}^\infty$ be a
zero-mean strong mixing process, where
$$
 \sum_{k=1}^{\infty}\alpha_X(k)^{\delta/(2+\delta)}<\infty, \ \ \mbox{and} \ \
  \sup_i E|X_i|^{2+\delta}=M<\infty, \  \ \mbox{for some} \ \delta\in (0, \infty).
$$
Then $\frac{E(S_n)^2}{n}$ is convergent.\\}

\noindent \textbf{Proof of Lemma 2:} By lemma 1,
$$
E(S_n)^2= \sum_{t=1}^{n}\sum_{s=1}^{n}E(X_tX_s) \leq
\sum_{t=1}^{n}\sum_{s=1}^{n}
[E|X_t|^p]^{1/p}[E|X_s|^q]^{1/q}[\alpha_X(s-t)]^{1-1/p-1/q}.
$$
Let $p=q=2+\delta$, then
$$
\frac{1}{n}E(S_n)^2 \leq
\frac{1}{n}M^{\frac{2}{2+\delta}}\sum_{i=-n}^{n}n[\alpha_X(i)]^{\delta/
2+ \delta} \leq M \sum_{i=-\infty}^{\infty}[\alpha_X(i)]^{\delta/
2+ \delta} \leq \infty.
$$
An alternative proof for lemma 2 can be found in Rio (1993).\\

\noindent
\textbf{Proof of Theorem 1:}
Let $$
 g_n(t):=\frac{[nt]}{\sqrt{n}}C^{-1/2}
(\hat{\gamma}_{[nt]}(0)-{\gamma}(0),
\hat{\gamma}_{[nt]}(1)-{\gamma}(1),...,
\hat{\gamma}_{[nt]}(L)-{\gamma}(L))^T,  \eqno(2.6)
$$
where $ \gamma(.)$ is auto-covariance function of
$\{X_t\}$.
 By   (2.4) and (2.6), it is immediate that
$$
g^*_n(t)= g_n(t) - \frac{[nt]}{n}g_n(1).
$$
Any Brownian bridge $W^0(t)$ has the same distribution as $W(t)-tW(1)$,
where $W(t)$ is a standard Brownian motion (Billingsley, 1999).\\
Therefore if $g_n(t) \Rightarrow (W_0(t), W_1(t), ..., W_L(t))^T$,
where $W_i(t)$ are independent Brownian motions for $0\leq i\leq L$,
then $g^*_n(t) \Rightarrow (W^0_0(t), W^0_1(t), ...,
W^0_L(t))^T$. So it is enough to show that $g_n(t)$ converges
to a vector of Brownian motions. \\
As
$$ g_n(t)= \frac{[nt]}{\sqrt{n}}C^{-1/2}
   \left( \begin{array}{c}
    \hat{\gamma}_{[nt]}(0)- \gamma(0)\\
    \vdots \\
    \hat{\gamma}_{[nt]}(L)- \gamma(L)\\
    \end{array} \right)
    =\frac{C^{-1/2}}{\sqrt{n}}
   \left( \begin{array}{c}
    \sum_{i=1}^{[nt]}X_{i}X_{i}- [nt]\gamma(0)\\
    \vdots \\
   \sum_{i=1}^{[nt]-L}X_{i}X_{i+L}- [nt]\gamma(L) \\
    \end{array} \right),
$$
so by assuming
$$
Y_{m,t}=X_tX_{t+m}- \gamma(m), \ \ \ \mbox{for} \ \ 0\leq m \leq L \eqno(2.7)
$$
we have that
$$
g_n(t)= \frac{C^{-1/2}}{\sqrt{n}}
   \left( \begin{array}{c}
    \sum_{i=1}^{[nt]}Y_{m,i}\\
    \vdots \\
   \sum_{i=1}^{[nt]-L}Y_{m,i}\\
    \end{array} \right)
    + \frac{C^{-1/2}}{\sqrt{n}}
    \left( \begin{array}{c}
    0\\
    \gamma(1)\\
    2\gamma(2)\\
    \vdots \\
    L\gamma(L)\\
    \end{array} \right). \eqno(2.8)
$$
By (2.7),
  $$
{\mathbf{M}}_{-\infty}^{n}(Y)=\sigma\{Y_{m,t}; t<n\}=\sigma\{X_tX_{t+m}; t<n\}\subseteq {\mathbf{M}}_{-\infty}^{n+m}(X),
  $$  and
$$
{\mathbf{M}}_{n+j}^{\infty}(Y)=\sigma\{Y_{m,t}; t>n+j\}=\sigma\{X_tX_{t+m}; t>n+j\}\subseteq {\mathbf{M}}_{n+j}^{\infty}(X).
$$
So relation (1.1) implies that, $\alpha_{Y_m}(j) \leq \alpha_X(j-m)$.
Hence $ \{ Y_{m,t}\}$ form a zero mean, strong
mixing process, where by assumption (ii),
$$
 \sum_{k=1}^{\infty}\alpha_{Y_m}(k)^{\delta/(2+\delta)}<\infty.  \
 \ \eqno(2.9)
 $$
Also by (2.7) and assumption (i),
$$
\sup_iE|Y_{m,i}|^{2+\delta}<\infty,  \ \ \eqno(2.10)
$$
for some $\delta\in(0,\infty)$. Let $S_{m,n}=Y_{m,1}+ ... +Y_{m,n} $.
Using (2.9) and (2.10), lemma 2
implies that $\frac{Var(S_{m,n})}{n}\rightarrow \sigma_m^2$,
for some $\sigma_m^2< \infty$.\\
If $\sigma_m> 0$, then (2.9), (2.10), and functional
central limit theorem, introduced by Herrndorf(1985), assert that
$$
 \frac{S_{m,[nt]}}{\sigma_m \sqrt{n}} \Rightarrow W(t), \ \
 \ \mbox{for} \ \ 0\leq m\leq L. \eqno(2.11)
$$
For $0\leq h, k\leq L$,
$$
\mbox{cov}(\frac{1}{\sqrt{n}}\sum_{i=1}^{[nt]-h}Y_{h,i},
\frac{1}{\sqrt{n}}\sum_{i=1}^{[nt]-k}Y_{k,i})=
 \mbox{cov}(\frac{1}{\sqrt{n}}\sum_{i=1}^{[nt]-h}Y_{h,i}+\frac{1}{\sqrt{n}}h\gamma(h),
 \frac{1}{\sqrt{n}}\sum_{i=1}^{[nt]-k}Y_{k,i}+\frac{1}{\sqrt{n}}k\gamma(k))
$$
$$
=\frac{[nt]^2}{n}\mbox{cov}(\hat{\gamma}_{[nt]}(h)- \gamma(h), \hat{\gamma}_{[nt]}(k)- \gamma(k))=
\frac{[nt]}{n}[nt]\mbox{cov}(\hat{\gamma}_{[nt]}(h), \hat{\gamma}_{[nt]}(k)).
$$

By Bartlett's formula (2.3),
$\lim_{n\rightarrow \infty}[nt]\mbox{cov}(\hat{\gamma}_{[nt]}(h), \hat{\gamma}_{[nt]}(k))=c_{hk}$.
So
$$
\lim_{n\rightarrow \infty}\mbox{cov}(\frac{1}{\sqrt{n}}\sum_{i=1}^{[nt]-h}Y_{h,i},
\frac{1}{\sqrt{n}}\sum_{i=1}^{[nt]-k}Y_{k,i})= t\,c_{hk}.  \eqno(2.12)
$$
 Hence by (2.11) and (2.12),
$$
\frac{C^{-1/2}}{\sqrt{n}}
   \left( \begin{array}{c}
    \sum_{i=1}^{[nt]}Y_{m,i}\\
    \vdots \\
   \sum_{i=1}^{[nt]-L}Y_{m,i}\\
    \end{array} \right) \Rightarrow \left( \begin{array}{c}
    W_0(t)\\
    \vdots \\
    W_L(t)\\
    \end{array} \right) \eqno(2.13)
$$
where $W_j(t), j=0, ..., L$, are independent Brownian motions.\\
The second part on right hand of (2.8) tends to zero as
 $n \rightarrow \infty$, so (2.13) implies that
$$
g_n(t)\Rightarrow
   \left( \begin{array}{c}
    W_0(t)\\
    \vdots \\
    W_L(t)\\
    \end{array} \right), \\ \ \  \
$$
where $W_j, j=0, ..., L$, are independent Brownian motions.$\Box$

\subsection{CUSUM test statistic}
\hspace{.2in} Using theorem 1 and (2.4), a CUSUM test statistic is
 constructed as:
$$
T_n := \max_{L \leq k < n}g^*_n(\frac{k}{n})^T .
g^*_n(\frac{k}{n}). \eqno(2.14)
$$
By continuous mapping theorem
$$
T_n  \Rightarrow \sup_{0\leq t< 1}\sum_{j=0}^{L}(W^0_j(t))^2.
\eqno(2.15)
$$
For detecting changes in time series $\{ X_t\}$, under the assumptions
of theorem 1, the following test is proposed for testing hypothesis

$H_0$: no change occur in the auto-covariance function of $X_1,
..., X_n$

$H_1$: there is a $1<k<n$  such that auto-covariance function of
$X_1, ..., X_k$ is different from auto-covariance function of
$X_{k+1}, ..., X_n.$

The strategy of this test is to reject $H_0$ when $T_n$ is large.\par
\hspace{.2in}By (2.14) and (2.15), the critical region of the test
at significant level $\alpha$ is $\{T_n \geq
c_{\alpha}\}$, where $c_{\alpha}$ is the $(1-\alpha)$-quantile
point of the distribution of $\sup_{0\leq t\leq 1}\sum_{j=0}^{L}(W^0_j(t))^2$.
The critical values can be found in Kiefer(1959) and Lee et al.(2003).

\textbf{Example:}
Let $\{X_t\}$ be an MA(1) process as:
$$
X_t=Z_t + \theta Z_{t-1},
$$
where $\{Z_t\}$ is an iid normal sequence with mean zero
and variance $\sigma^2$. In linear processes where $E(Z_t^4)=\eta\sigma^4$,
 Bartlett's formula has explicit form as:
$$
c_{i,j}= \Sigma_{l=-\infty}^{\infty}\{\gamma(l)\gamma(l-i+j) + \gamma(l+j)\gamma(l-i)\}
+(\eta-3)\gamma(i)\gamma(j),
$$
where $\gamma(l)$ is the corresponding autocovariance function at lag $l$
of $\{X_t\}$, see Brockwell and Davis(1991).\\
If the noise is Gaussian, $\eta=3$, so
$$
c_{i,j}=\gamma(1)\gamma(j-i-1)+\gamma(1+i)\gamma(j-1)+\gamma(0)\gamma(j-i)
+\gamma(i)\gamma(j)+\gamma(1)\gamma(j-i+1)
$$
$$
+\gamma(1+j)\gamma(1-i).
$$
\input{epsf}
\epsfxsize=2.5in \epsfysize=2in
\begin{figure}
\caption{\scriptsize behavior of $T_n$ for different values of parameter $\theta$ in MA(1) model}
\centerline{\epsfxsize=4in \epsfysize=2in \epsffile{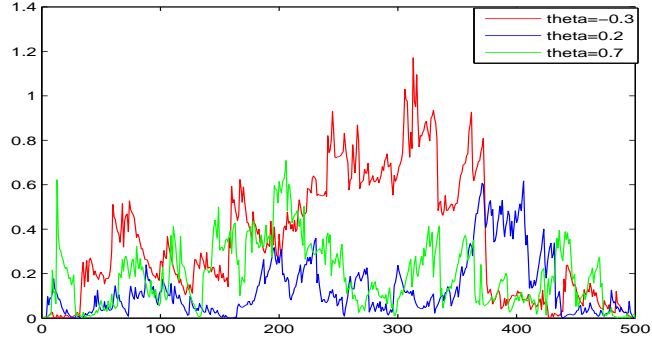}}
\end{figure}
Let $L=1$, therefore corresponding covariance matrix
 $C=[c_{i,j}]_{i,j=1}^{2}$ can be written as
$$
C=\left(
  \begin{array}{cc}
   2\gamma^2(0)+4\gamma^2(1) & 4\gamma(0)\gamma(1) \\
   4\gamma(0)\gamma(1)       & \gamma^2(0)+3\gamma^2(1) \\
  \end{array}
\right)= \left(
  \begin{array}{cc}
   2(1+4\theta^2+\theta^4) & 4\theta(1+\theta^2) \\
   4\theta(1+\theta^2)     & (1+5\theta^2+\theta^4) \\
  \end{array}
\right)\sigma^2.
$$
By (2.4), for $t=\frac{k}{n}$ we have that
$$
g^*_n(k/n):=\frac{k}{\sqrt{n}}C^{-1/2}(
\hat{\gamma}_{k}(0)-\hat{\gamma}_n(0),\,
\hat{\gamma}_{k}(1)-\hat{\gamma}_n(1))^T,
$$
and  by (2.15)
$$
T_n=\max_{L \leq k < n}\{\frac{k^2}{n}(\hat{\gamma}_{k}(0)-\hat{\gamma}_n(0),\,
\hat{\gamma}_{k}(1)-\hat{\gamma}_n(1))C^{-1}(\hat{\gamma}_{k}(0)-\hat{\gamma}_n(0),\,
\hat{\gamma}_{k}(1)-\hat{\gamma}_n(1))^T\}.
$$
\hspace{.2in}Figure 1 shows the behavior of test statistic $T_n$
for different values of parameter $\theta$ in an MA(1) process without
change point. For $L=1$, at significant level $\alpha=5\%$, the
critical value is $c_{\alpha}=2.408$.

\subsection{Consistent estimation of covariance matrix}
\hspace{.2in}In this section we evaluate asymptotic behavior of the covariance
function of estimators $\hat{\gamma}_n(i)$, defined by (2.3).
Bartlett(1946) derived an explicit formula for the asymptotic
behavior of covariance function of sample autocovariances,
when there exist a linear model for the data(Priestley, 1981).
We present a consistent estimator $\hat{C}$ for the case that
there is no model for data, or we have nonlinear process.
Using (2.3), Bartlett's estimator can be written as:
$$
c_{hk}=\lim_{n\rightarrow \infty}\theta_n(h, k)=\lim_{n\rightarrow \infty}n\,\mbox{cov}(\hat{\gamma}_n(h),\hat{\gamma}_n(k)), \ \
0\leq h \leq k <n.
$$
\hspace{.2in}For stationary process $\{ X_t\}$, let
$$
 \tilde{\gamma}_n(h)=\frac{1}{n}\sum_{t=1}^{n}X_tX_{t+h}, \ \ 0\leq
 h<n, \eqno(2.16)
$$
and $\tilde{\theta}_n(h, k)
=n\,\mbox{cov}(\tilde{\gamma}_n(h),\tilde{\gamma}_n(k))$.
By theorem 3 it is shown that $\tilde{\theta}_n(h, k)$ has
the same asymptotic behavior  as $\theta_n(h, k)$.
For the evaluation of $\tilde{\theta}_n(h, k)$,
one can easily verify that
$$
\tilde{\theta}_n(h, k)
=n\,\mbox{cov}(\tilde{\gamma}_n(h),\tilde{\gamma}_n(k))=
n\{E(\tilde{\gamma}_n(h).\tilde{\gamma}_n(k))-\gamma(h).\gamma(k)\}
$$
$$
=nE\frac{1}{n^2}\sum_{t=1}^{n}\sum_{s=1}^{n}\{X_tX_{t+h}X_sX_{s+k}
-\gamma(h)\gamma(k)\}$$
$$
=\frac{1}{n}\sum_{t=1}^{n}\sum_{s=1}^{n}
E\{X_tX_{t+h}X_sX_{s+k}-\gamma(h)\gamma(k)\}=\frac{1}{n}\sum_{l=0}^{n-1}
\sum_{t=1}^{n-l}\{ E(X_tX_{t+h}X_{t+l}X_{t+l+k})
$$
$$
-\gamma(h)\gamma(k)\} +
\frac{1}{n}\sum_{l=-n+1}^{-1}\sum_{t=-l+1}^{n}\{
E(X_tX_{t+h}X_{t+l}X_{t+l+k})-\gamma(h)\gamma(k)\}.
$$
By replacing $l$ with $-l$ and $t-l$ with $t$ in last summation we have
$$
\tilde{\theta}_n(h, k)=\frac{1}{n}\sum_{l=0}^{n-1}\sigma_{h,k}(l),
\eqno(2.17)$$
where
$$
\sigma_{h,k}(0)=\sum_{t=1}^{n}\{E(X_tX_{t+h}X_{t}X_{t+k})
-\gamma(h)\gamma(k)\}, \eqno(2.18)
$$
and for $1\leq l \leq n-1$
$$
\sigma_{h,k}(l)=\sum_{t=1}^{n-l}\{E(Y_{1t}^{l})+
E(Y_{2t}^{l})-2\gamma(h)\gamma(k)\}, \eqno(2.19)
$$
in which  $Y_{1t}^{l}= X_tX_{t+h}X_{t+l}X_{t+l+k}$ and
$Y_{2t}^{l}= X_{t+l}X_{t+l+h}X_{t}X_{t+k}$.\\

Let
$$\bar{\theta}_n(h, k)
=\frac{1}{n}\sum_{l=0}^{h_n}\bar{\sigma}_{h,k}(l),
 \eqno(2.20) $$
where
$$
\bar{\sigma}_{h,k}(0)= \sum_{t=1}^{n}\{\frac{1}{n}\sum_{i=1}^{n}\{Y_{1i}^{0}-
\tilde{\gamma}_n(h)\tilde{\gamma}_n(k)\}, \quad \quad \quad \eqno(2.21)
$$
$$
\bar{\sigma}_{h,k}(l)= \sum_{t=1}^{n-
l}\{\frac{1}{n}\sum_{i=1}^{n}\{Y_{1i}^{l} + Y_{2i}^{l}\}-
2\tilde{\gamma}_n(h)\tilde{\gamma}_n(k)\},  \quad \quad  1\leq l <n \eqno(2.22)
$$
and $\{h_n\}$ is a sequence of positive integers that
$$
h_n=O(n^{\beta}) \ \ \;  \mbox{for} \ \, \mbox{some} \; \ \beta
\in (0, 1/2). \eqno(2.23)
$$
Now we have the following result.\\
\emph{ \noindent
\textbf{Theorem 2}: Under the assumptions of theorem 1, if  \
$ \sup_tE|X_t|^{8+\delta}< \infty$ for some $\delta>0$  then
$$
\parallel \bar{\theta}_n(h, k)- \tilde{\theta}_n(h, k)
\parallel_2 \ \rightarrow 0 \ \ as \ n\rightarrow \infty.
$$
}
\noindent
\textbf{Proof :} The proof is organized in three
steps:
\par
\hspace{.2in}\textmd{Step 1:}
$\lim_{n\rightarrow \infty}\frac{1}{n}\sum_{1\leq l< n}|\sigma_{h,k}(l)|<\infty$.

By lemma 1 and (2.19),
$$|\sigma_{h,k}(l)|\leq |\sum_{t=1}^{n-l}\{E(X_tX_{t+h}X_{t+l}X_{t+l+k})
- \gamma(h)\gamma(k)\}| +
|\sum_{t=1}^{n-l}\{E(X_tX_{t+h}X_{t+l}X_{t+l+k}) -
\gamma(h)\gamma(k)\}|
$$
$$
\leq
\sum_{t=1}^{n-l}[E|X_tX_{t+h}|^p]^{1/p}[E|X_{t+l}X_{t+l+k}|^q]^{1/q}[\alpha(l-h)]^{1-1/p-1/q}
+$$
$$
\sum_{t=1}^{n-l}[E|X_{t+l}X_{t+l+h}|^p]^{1/p}[E|X_{t}X_{t+k}|^q]^{1/q}[\alpha(l+h)]^{1-1/p-1/q}.
$$
So for $p=q=2+\delta$, by assumption (ii) of theorem 1,
$$
|\sigma_{h,k}(l)|\leq M(n-l)\{
\alpha_X(l-h)^{\frac{\delta}{2+\delta}}+
\alpha_X(l+h)^{\frac{\delta}{2+\delta}}\},
$$
and
$$
\lim_{n\rightarrow \infty}\frac{1}{n}\sum_{1\leq l< n}|\sigma_{h,k}(l)|<\infty.
\eqno(2.24)
$$
\par
\hspace{.2in}\textmd{Step 2:}
$\|\bar{\sigma}(l) - \sigma(l)\|_2$ is of order $O(n^{1/2})$.\\
For $1 \leq l < n $, by (2.19) and (2.22)
$$
\|\sigma_{h,k}(l) -\bar{\sigma}_{h,k}(l) \|_2=
\|\sum_{t=1}^{n-l}\{EY_{1t}^l + EY_{2t}^l - 2\gamma(h)\gamma(k)\}
- \sum_{t=1}^{n-l}\{\frac{1}{n}\sum_{i=1}^{n}Y_{1i}^l +
\frac{1}{n}\sum_{i=1}^{n}Y_{1i}^l-2\tilde{\gamma}_n(h)\tilde{\gamma}_n(k)\}\|_2
$$
$$
=\|\sum_{t=1}^{n-l}\{EY_{1t}^l + EY_{2t}^l - 2\gamma(h)\gamma(k) -
\frac{1}{n}\sum_{i=1}^{n}Y_{1i}^l +
\frac{1}{n}\sum_{i=1}^{n}Y_{1i}^l -
2\tilde{\gamma}_n(h)\tilde{\gamma}_n(k)\}\|_2
$$
$$
 \leq \sum_{t=1}^{n-l}\{a_1 + a_2 + 2a_3\},
$$
\\
where $ a_1= \|\frac{1}{n}\sum_{i=1}^{n}\{Y_{1i}^l -
EY_{1t}^l\}\|_2 $, $ a_2= \|\frac{1}{n}\sum_{i=1}^{n}\{Y_{2i}^l -
EY_{2t}^l\}\|_2 $ and $ a_3= \|\tilde{\gamma}_n(h)\tilde{\gamma}_n(k)-
\gamma(h)\gamma(k)\|_2$. As by the assumption of the theorem
$\sup_tE|Y_t|^{2+\delta/4}<\infty$, so by lemma 2
$$
a_1^2= \frac{1}{n^2}E( \sum_{i=1}^{n}\{Y_{1i}^l - EY_{1t}^l\})^2=
O(n^{-1}),
$$
and $a_1=O(n^{-1/2})$. Similarly $a_2=O(n^{-1/2})$. Also
$$
a_3=\|\tilde{\gamma}_n(h)\tilde{\gamma}_n(k)- \gamma(h)\gamma(k)\|_2=
\|\tilde{\gamma}_n(h)\tilde{\gamma}_n(k)- \tilde{\gamma}_n(h)\gamma(k) +
\tilde{\gamma}_n(h)\gamma(k)- \gamma(h)\gamma(k)\|_2
$$
$$
\leq \|\tilde{\gamma}_n(h)\|_2 \|\tilde{\gamma}_n(k)- \gamma(k)\|_2 +
\|\tilde{\gamma}_n(h)- \gamma(h)\|_2 |\gamma(k)|. \eqno(2.25)
$$
Also by lemma 2,
$$
\|\tilde{\gamma}_n(h)\|_2^2=
\|\frac{1}{n}\sum_{t=1}^{n}X_tX_{t+h}\|_2^2=
\frac{1}{n}E(\sum_{t=1}^{n}X_tX_{t+h})^2 < \infty, \eqno(2.26)
$$
and
$$
\|\tilde{\gamma}_n(k)- \gamma(k)\|_2^2=
\|\frac{1}{n}\sum_{t=1}^{n}X_tX_{t+h}- \gamma(k)\|_2^2 \leq
\frac{1}{n^2}E(\sum_{t=1}^{n}X_tX_{t+h}- \gamma(k))^2= O(n^{-1}).
\eqno(2.27)
$$
Thus by relations (2.25), (2.26) and (2.27) we have that $a_3= O(n^{-1/2})$. \\
By similar method, one can easily verify that
 $\|\bar{\sigma}_{h, k}(0) - \sigma_{h, k}(0)\|_2=O(n^{1/2})$. Therefore,
$$
\|\bar{\sigma}_{h, k}(l) - \sigma_{h, k}(l)\|_2 \leq O(n^{1/2}).\eqno(2.28)
$$
\par
\hspace{.2in}\textmd{Step 3:} Steps 1 and 2 are applied to prove the main result.\\
By Minkowski inequality we have
$$
 \| \tilde{\theta}_n(h, k) -\bar{\theta}_n(h, k) \|_2 =
  \|\frac{1}{n}\sum_{l=1}^{n}\sigma_{h, k}(l)-
   \frac{1}{n}\sum_{l=1}^{h_n}\bar{\sigma}_{h, k}(l)\|_2\leq
$$
$$
\frac{1}{n}\sum_{l=1}^{h_n}
 \|\sigma_{h, k}(l)- \bar{\sigma}_{h, k}(l)\|_2 \ +
 \frac{1}{n}\sum_{l=h_n}^{n}|\sigma_{h, k}(l)|. \eqno(2.29)
$$
So by (2.23) and (2.28) \\
$$
\frac{1}{n}\sum_{l=1}^{h_n}
 \|\sigma_{h, k}(l)- \bar{\sigma}_{h, k}(l)\|_2 \
  \leq \frac{1}{n}h_nO(n^{1/2})\rightarrow 0, \ \  \ as \ \
n\rightarrow\infty. \eqno(2.30)
$$
As $h_n\rightarrow \infty$, (2.24) implies that
$$
 \sum_{h_n< l< n}\frac{1}{n}|\sigma_{h,k}(l)| \rightarrow 0, \ \
as \ n\rightarrow \infty. \eqno(2.31)
$$
Finally by (2.30) and (2.31), we arrive at the assertion of the theorem.
$\Box$

\emph{\noindent
\textbf{Theorem 3}: Under the assumptions of theorem 2,
we have that
 $\| \theta_n(h, k) - \bar{\theta}_n(h, k)\|_2 \rightarrow 0, \ as\ n \rightarrow\infty
 $.\\}
\noindent
\textbf{Proof }: As
$$
\| \theta_n(h, k) - \bar{\theta}_n(h, k)\|_2\leq \| \theta_n(h,
k) - \tilde{\theta}_n(h, k)\|_2 + \| \tilde{\theta}_n(h, k) -
\bar{\theta}_n(h, k)\|_2,
$$
so by theorem 2 the second part on the right tends to zero, and for the first part,
by (2.2) and (2.16)
$$\| \theta_n(h, k) - \tilde{\theta}_n(h, k)\|_2= n\|\mbox{cov}(\hat{\gamma}_n(h),\hat{\gamma}_n(k))-
\mbox{cov}(\tilde{\gamma}_n(h),\tilde{\gamma}_n(k))\|_2
$$
$$ = \frac{1}{n}\| \sum_{t=1}^{n-h}
\sum_{s=1}^{n-k}E\{X_tX_{t+h}X_{s}X_{s+k}- \gamma(h)\gamma(k)\} -
\sum_{t=1}^{n} \sum_{s=1}^{n}E\{X_tX_{t+h}X_{s}X_{s+k}-
\gamma(h)\gamma(k)\}\|_2
$$
$$
=\frac{1}{n}\sum_{t=n-h+1}^{n}\sum_{s=n-k+1}^{n}\|\{X_tX_{t+h}X_{s}X_{s+k}-
\gamma(h)\gamma(k)\}\|_2
 \rightarrow 0,\  \mbox{as} \ n\rightarrow \infty.\Box
$$
\\
\noindent
\textbf{Corollary 1}: Under the assumptions of theorem 2,
by choosing $\hat{C}=[\hat{c}_{hk}]_{L+1\times L+1}$, where
  $\hat{c}_{hk}= \bar{\theta}_n(h,k)$ defined by (2.20),
$\hat{C}$ is a consistent estimator of covariance matrix $C$ in relation
(2.4).

\noindent
\textbf{Corollary 2}: Under the assumptions of theorem 2,
by (2.4), (2.14), (2.15), and corollary 1, we have
$$
 \hat{T}_n= \max_{1\leq k\leq n}\hat{g}^*_n(k/n)^T . \hat{g}^*_n(k/n) \Rightarrow \sup_{0\leq s \leq 1}\sum_{j=0}^{L}(W^0_j(s))^2,
 \eqno(2.32)
$$
where
$$
\hat{g}^*_n(s):=\frac{[ns]}{\sqrt{n}}\hat{C}^{-1/2}
(\hat{\gamma}_{[ns]}(0)-\hat{\gamma}_n(0),
\hat{\gamma}_{[ns]}(1)-\hat{\gamma}_n(1),...,
\hat{\gamma}_{[ns]}(L)-\hat{\gamma}_n(L))^T, \eqno(2.33)
 $$
in which $\hat{C}=[\bar{\theta}_n(h,k)]$, and $\bar{\theta}_n(h,k)$
is defined by (2.20).

\section{Simulation results}
\hspace{.2in}In this section, we investigate the performance of the proposed
 test statistic $\hat{T}_n$, by  a simulation study.
As this test statistic is to be applied for linear and nonlinear models,
 we consider simulations of such classes of time series.
 Test statistics are evaluated by using
  relations (2.32) and (2.33) with $L=1$ and relation (2.20) with $h_n =n^{0.3}$. \par
\hspace{.2in}For creating change in the covariance structure of time series,
the parameters are changed at the midpoint of the series.
 Empirical powers are evaluated, and for the case that there is
 no change in data, probability of type $I$ error is reported.\par
\hspace{.2in} The critical value of the test statistic at level   $\alpha= 0.05$
  is $c_{\alpha}=2.408$(Lee at el. 2003).
Simulations are repeated 1000 times, for  the following linear and nonlinear models, to
 evaluated empirical powers.\\
\textbf{Linear models:}
\begin{itemize}
    \item Model 1: \\
    ARMA(1,1): \ $X_t-\phi X_{t-1}= Z_t + \theta Z_{t-1}$,
    \item Model 2: \\
    MA(2): \ $X_t= Z_t + \theta_1 Z_{t-1}+ \theta_2 Z_{t-2}$,
\end{itemize}
 where $\{Z_t\}$ is a sequence of iid normal
random variables with mean zero and variance 1.\\
\begin{table}
\centering \caption{\scriptsize Empirical power of $\hat{T}_n$ for ARMA(1,1), where the initial parameter  $(\theta_0, \phi_0)=(0.1,0.2)$.}
\begin{tabular}{ c | cc cc cc cc }
\hline
&$\phi_1$&            &&            &&       &  \\
$\theta_1$&& 0.2 && 0.4 && 0.5 && 0.6     \\ \hline
 0.1  && $0.047^*$ && 0.408 && 0.761 && 0.935   \\
 0.3  && 0.295 && 0.874 && 0.968 && 0.989   \\
 0.5  && 0.734 && 0.977 && 0.994 && 0.998   \\
 0.7  && 0.935 && 0.995 && 0.999 && 1.000   \\ \hline
\end{tabular}
\end{table}
\hspace{.2in}
Table 1 reports results of the test for simulated data
from ARMA(1,1), where  250 samples are generated  with $(\theta_0,\phi_0)=(0.1, 0.2)$, and then 250 more samples
  for different values of $(\theta_1,\phi_1)$
 as reported in the table. In table 1, empirical powers are evaluated
 for cases where one or both parameters have changed.
In all cases high empirical powers shows ability of this test
statistic. The value pointed by $^*$ is type $I$ error, empirical level,
  which is slightly below $0.05$. \par
\input{epsf}
\epsfxsize=2.5in \epsfysize=2in
\begin{figure}
\caption{\scriptsize Empirical power of CSSM in comparison with Berkes et al.(2009)
in MA(2) model, $X_t=Z_t+ \theta Z_{t-1}+ \theta Z_{t-2}$.
 parameter $\theta$ changes from 0 to alternative $\theta$.}
\centerline{\epsfxsize=3in \epsfysize=2in \epsffile{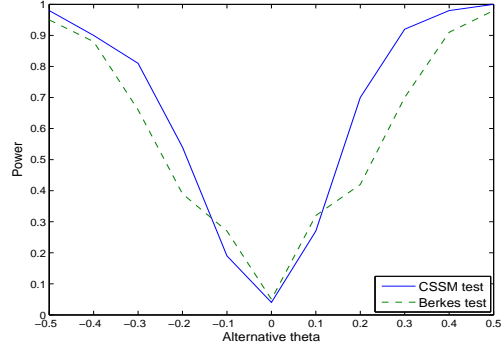}}
\end{figure}
\input{epsf}
\epsfxsize=2.5in \epsfysize=1.7in
\begin{figure}
\caption{\scriptsize Behavior of CSSM and Berkes statistic in simulated samples from MA(2) model,
without change point.}
\centerline{\epsffile{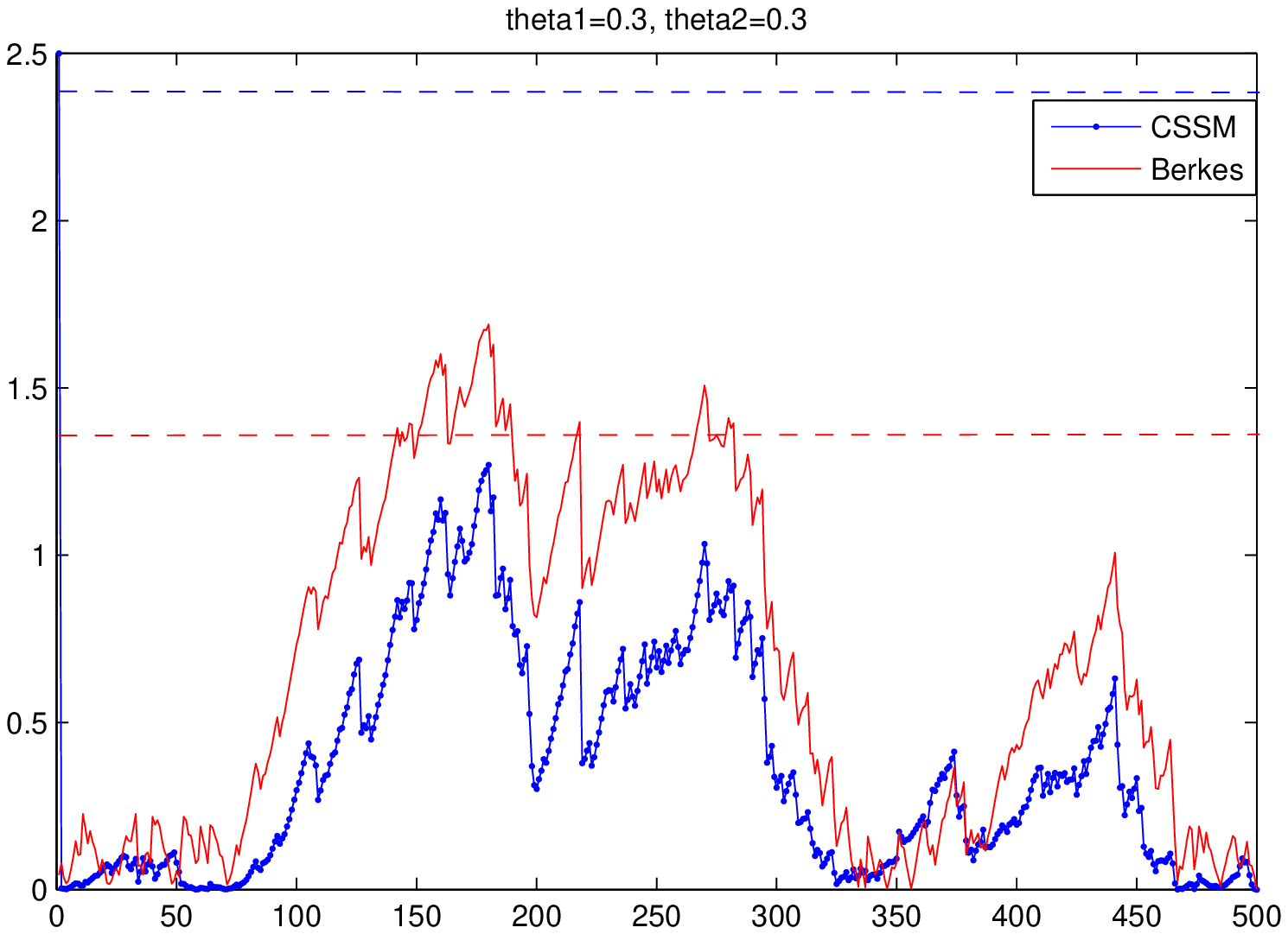} $\hspace{.2in}$ \epsfxsize=2.5in
\epsfysize=1.7in \epsffile{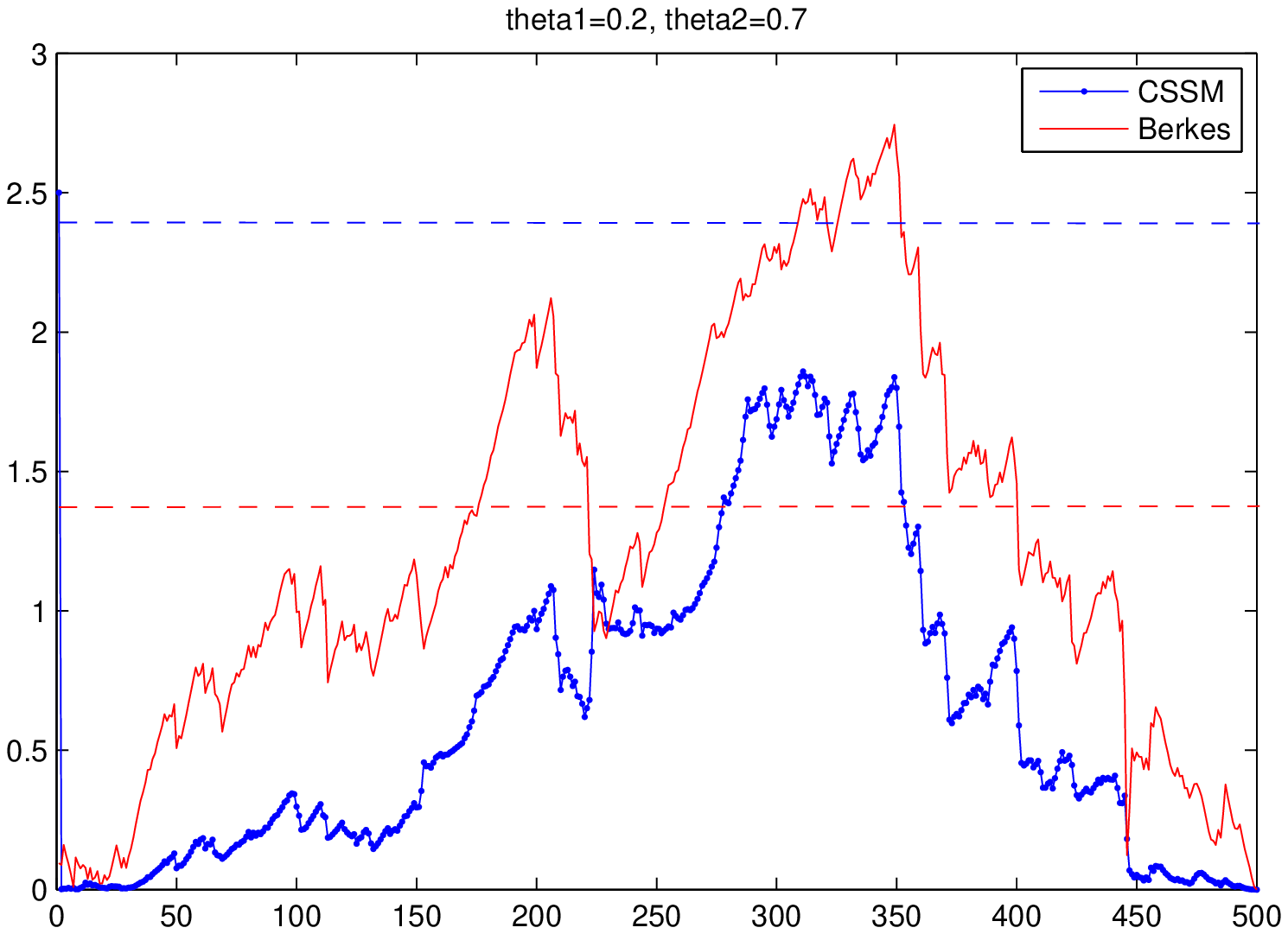}}
\end{figure}
\input{epsf}
\epsfxsize=2.5in \epsfysize=1.7in
\begin{figure}
\caption{\scriptsize Simulated samples from 2-dependent model with a change point in $k^*=500$, variance
of $Z_t$ change from 1 to 1.26(left), Corresponding test statistics(right).}
\centerline{\epsffile{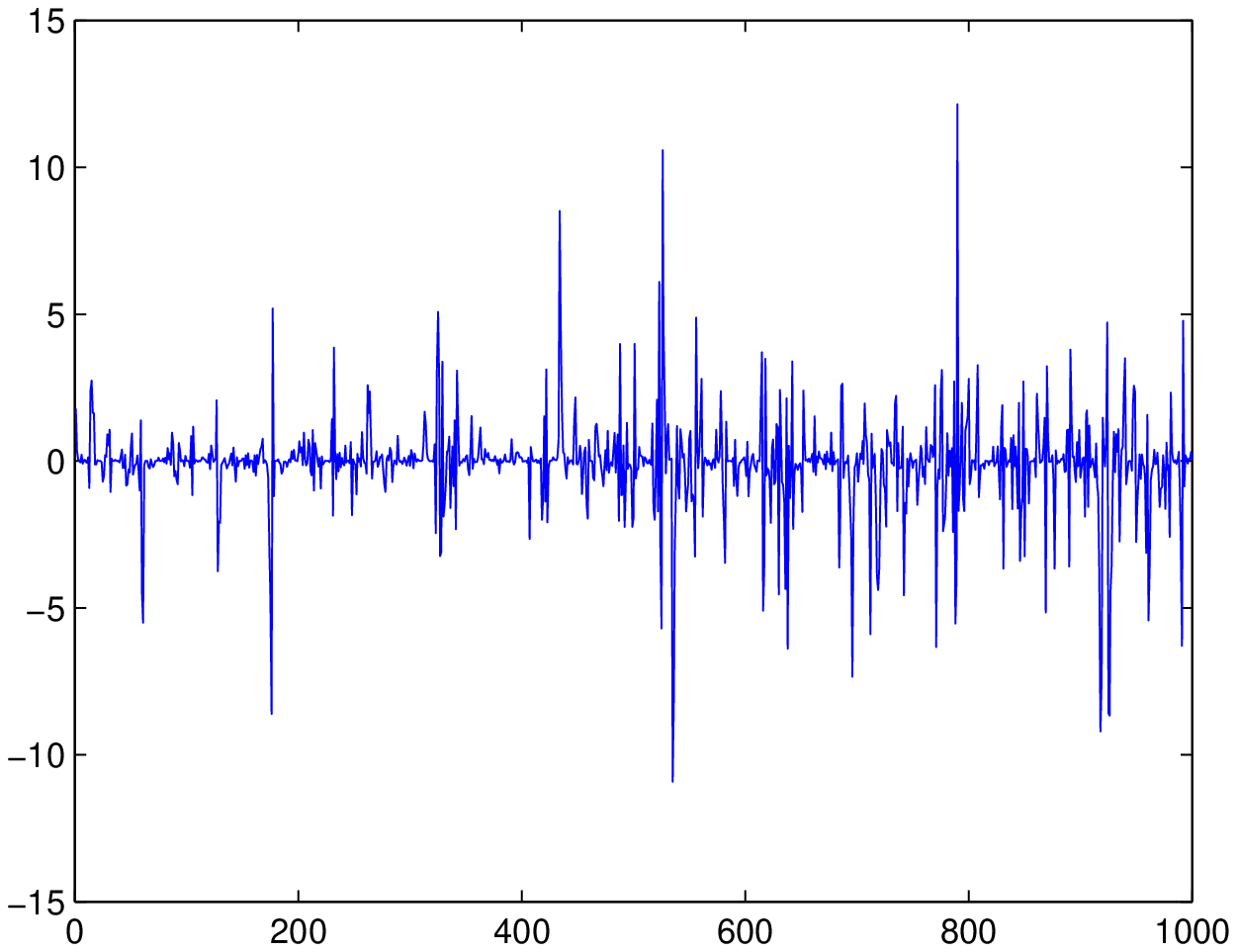} $\hspace{.2in}$ \epsfxsize=2.5in
\epsfysize=1.7in \epsffile{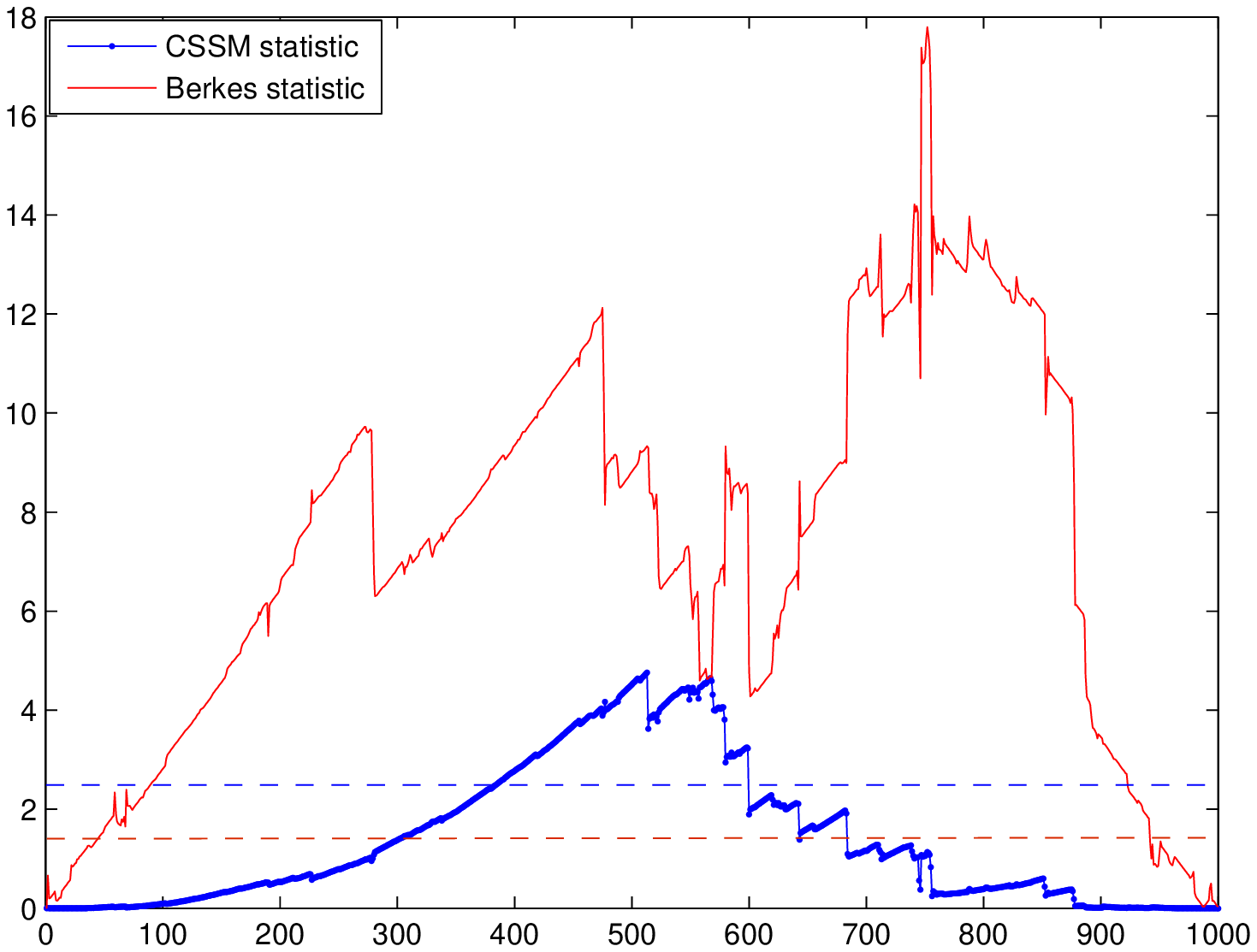}}
\end{figure}
\hspace{.2in} As a comparison of CSSM with the Berkes method for linear models,  mode 2
 with $\theta=\theta_1=\theta_2$ is considered and
 250 samples are generated with $\theta=0$ at first stage
and then 250 more samples for some alternative $\theta$,
and empirical powers are plotted in figure 2.
As empirical powers show, CSSM test has a
better performance, with smaller  type $I$ error.\par
 To visualize this, we follow a simulation of model 2, where there is no change in
 time series data. So we generate 500 samples from MA(2) with $(\theta_1,\theta_2)=(0.3,0.3)$
 once and again with $(\theta_1,\theta_2)=(0.2,0.7)$. Then we evaluate  CSSM and Berkes
 test statistics, and plot them with corresponding critical values, 2.408 and 1.36
 in figure 3. As figure 3 shows Berkes method indicate a change point by mistake
 where CSSM statistic is far beyond such miss detection.\\
\textbf{Nonlinear models:}\par
\hspace{.2in}As nonlinear models, we consider followings:
\begin{itemize}
  \item Model 3:\\
  2-dependent: \ $X_t=Z_tZ_{t-1}Z_{t-2}$
  \item Model 4:\\
  GARCH(1,1): \ $X_t=h_t.Z_t$\
  \hspace{.1in} where \  $h_t^2=\omega + \alpha X^2_{t-1}+ \beta h_{t-1}^2$
\end{itemize}
in which $\{Z_t\}$ is a sequence of iid normal
random variables with mean zero and variance one, $\sigma^2=1$. \par
\hspace{.2in}Figure 4(left) shows 1000 generated samples of a 2-dependent
process, model 3, with a change point at $k^*=500$, where variance
of $\{Z_t\}$ changes from 1 to 1.26.  Figure 4(right) shows the
behavior of CSSM and Berkes test statistics.
 Corresponding critical values  are 2.408 and
 1.36 respectively which are presented by horizontal lines.
 Figure 4(right) shows that the supremums  of both
statistics exceed corresponding  critical values and a change in process is detected, but CSSM statistics is more precise, as
CSSM  detects the change at t=512 and Berkes  statistic
detects it at t=750.  \par
\hspace{.2in}Table 2 reports the empirical power of the CSSM,
 for 2-dependent model, model 3.
 In table 2(a) change of the variance of $\{Z_t\}$, from
 $\sigma^2=1$ to alternative values is proposed.
In table 2(b) change of the mean of $\{Z_t\}$, from
$\mu_0=0$  to alternative values, has been considered.\\
\begin{table}
\centering \caption{ \scriptsize Empirical power in 2-dependent model.}
\begin{tabular}{ c c|c c }
\hline
Table 2(a) && Table 2(b) \\\hline
change in variance of $Z_t$&&  change in mean of $Z_t$ &   \\ \cline{1-2}\cline{3-4}
   $\sigma$& Power  &  $\mu$  &  Power           \\ \hline
      0.8    & 0.622  &   0.0   & 0.049          \\
      0.6    & 0.960  &   0.5   & 0.285          \\
      0.4    & 0.975  &   1.0   & 0.961          \\
      0.2    & 0.989  &   1.5   & 0.999          \\
\hline
\end{tabular}
\end{table}
\begin{table}
\centering \caption{ \scriptsize Empirical power of test in GARCH(1,1) model.}
\begin{tabular}{ c| c c c cccc c }
\hline
$(\omega, \alpha, \beta)$&n& 500   && 800  && 1000 \\ \hline
    no change            & & 0.034 && 0.035 && 0.032 \\
(0.8, 0.1, 0.2)          & & 0.528 && 0.748 && 0.894 \\
(0.8, 0.1, 0.5)          & & 0.735 && 0.931 && 0.967 \\
(0.8, 0.4, 0.2)          & & 0.974 && 0.999 && 1.000 \\
\hline
\end{tabular}
\end{table}
\hspace{.2in}Table 3 shows the empirical power of the CSSM statistics in a GARCH(1,1) model.
These simulations are  done for different values of $n$.
Parameters initial values are
$(\omega, \alpha, \beta)=(0.5, 0.1, 0.2)$, and empirical powers
for different alternatives are reported.
Simulations show that the powers has significant increase with sample size.

\section{Conclusion}
\hspace{.2in}In this article a nonparametric CUSUM test statistic is proposed
for detecting structural changes in strong mixing time series. Under
a sufficient condition this test
statistic converges in distribution to the supremum  of the sum of
independent standard Brownian bridges. This method covers a broad
class of linear and nonlinear time series such as ARMA and
GARCH models, m-dependent models and many others.
Beside the wide applications, simulation results
shows that our test statistic in comparison with Berkes et al.(2009)
has a better performance and  is more powerful.

\section{Acknowledgment}
The authors would like to thank the referee for careful reading of the paper and a number of valuable comments and suggestions that have improved the quality of the paper.

\end{document}